\documentclass{amsart}
\usepackage{amsfonts,amssymb,amscd,amsmath,enumerate,verbatim,newlfont,calc}

\newtheorem{theorem A}{Theorem A}
\newtheorem{theorem B}{Theorem B}

\begin{document}
\authors
\title{The Bogomolov multiplier of Lie algebras}
\author[Z. Araghi Rostami]{Zeinab Araghi Rostami}
\author[M. Parvizi]{Mohsen Parvizi}
\author[P. Niroomand]{Peyman Niroomand}
\address{Department of Pure Mathematics\\
Ferdowsi University of Mashhad, Mashhad, Iran}
\email{araghirostami@gmail.com, zeinabaraghirostami@stu.um.ac.ir}
\address{Department of Pure Mathematics\\
Ferdowsi University of Mashhad, Mashhad, Iran}
\email{parvizi@um.ac.ir}
\address{School of Mathematics and Computer Science\\
Damghan University, Damghan, Iran}
\email{niroomand@du.ac.ir, p$\_$niroomand@yahoo.com}
\maketitle
\begin{abstract}
In this paper, we extend the notion of the Bogomolov multipliers and the CP-extensions to Lie algebras. Then we compute the Bogomolov multipliers for Abelian, Heisenberg and nilpotent Lie algebras of class at most 6. Finally we compute the Bogomolov multipliers of some simple complex Lie algebras.
\end{abstract}
\section{\bf{Introduction}}
During the study of continuous transformation groups in the end of 19th century, Sophus Lie found new algebraic structures now known as \emph{Lie algebras}. This new structure played an important role in 19th and 20th centuries mathematical physics. (See [31,33], for more information). Also, today, more than a century after Lie's discovery, we have an extensive and important algebraic theory studying objects like Lie algebras, Lie groups, Root systems, Weyl groups, Linear algebraic groups, etc; which is named \emph{Lie theory} and the some researchs show its emphasis in modern mathematics. (See [5,31], for more information). Furthermore, mathematicians discovered that every Lie algebra could be associated to a continuous or Lie group. For example, Lazard introduced a correspondence between some groups and some Lie algebras. (See [21], for more information). So, theories of groups and Lie algebras are structurally similar, and many concepts related to groups, are defined  analogously to Lie algebras. In this paper we want to define the Bogomolov multipliers for Lie algebras. This concept is known for groups. The Bogomolov multiplier is a group-theorical invariant introduced as an obstruction to the rationality problem in algebraic geometry. Let $K$ be a field, $G$ a finite group and $V$ a faithful representation of $G$ over $K$. Then there is natural action of $G$ upon the field of rational functions $K(V)$. The rationality problem (also known as Noether's problem) asks whether the field of $G$-invariant functions ${K(V)}^G$ is rational (purely transcendental) over $K$? A question related to the above mentioned is whether there exist indpendent variables $x_1,...,x_r$ such that ${K(V)}^{G}(x_1,...,x_r)$ becomes a pure transcendental extension of $K$? Saltman in [29] found examples of groups of order $p^9$ with negative answer to the Noether's problem, even when taking $K=\Bbb{C}$. His main method was the application of the unramified cohomology group ${H_{nr}^{2}}({\Bbb{C}(V)}^{G},{\Bbb{Q}}/{\Bbb{Z}})$ as an obstruction. Bogomolov in [4] proved that it is canonically isomorphic to
$$B_0(G)={\bigcap}\ {\ker} \{ res_{G}^{A} : H^2(G,{\Bbb{Q}}/{\Bbb{Z}}) \rightarrow H^2(A,{\Bbb{Q}}/{\Bbb{Z}}) \},$$
where $A$ is an abelian subgroup of $G$.
The group $B_0(G)$ is a subgroup of the Schur multiplier $\mathcal{M}(G)=H^2(G,{\Bbb{Q}}/{\Bbb{Z}})$ and Kunyavskii in [20] named it the \emph{Bogomolov multiplier} of $G$. Thus vanishing the Bogomolov multiplier leads to positive answer to Noether’s problem. But it's not always easy to calculate Bogomolov multipliers of groups.  Recently, Moravec in [26] introduced an equivalent definition of the Bogomolov multiplier. In this sense, he used a notion of the nonabelian exterior square of a group $G$ ($G\wedge G$) to obtain a new description of Bogomolov multiplier. He showed that if $G$ is a finite group, then $B_0(G)$ is non-canonically isomorphic to $\text{Hom} (\tilde{B_0}(G),{\Bbb{Q}}/{\Bbb{Z}})$, where the group $\tilde{B_0}(G)$ can be described as a section of the nonabelian exterior square $G\wedge G$ of the group $G$. Also, he proved that $\tilde{B_0}(G)=\mathcal{M}(G)/\mathcal{M}_0(G)$, such that the Schur multiplier $\mathcal{M}(G)$ or the same $H^2(G,{\Bbb{Q}}/{\Bbb{Z}})$ interpreted as the kernel of the commutator homomorphism $G\wedge G \rightarrow [G,G]$ given by $x\wedge y \rightarrow  [x,y]$, and $\mathcal{M}_0(G)$ is a subgroup of $\mathcal{M}(G)$ defined as $\mathcal{M}_0(G)=<x\wedge y \ | \  [x,y]=0 , \  x,y\in G>$. Thus, in finite case, $\tilde{B_0}(G)$ is non-canonically isomorphic to $B_0(G)$.
 With this definition and similar to the Schur multiplier, the Bogomolov multiplier can be explicated as a measure of the extent to which relations among commutators in a group fail to be consequences of universal relation. Furthermore, Moravec's method relates Bogomolov multiplier to commuting probability of a group and shows that the Bogomolov multiplier plays an important role in commutativity preserving central extensions of groups, that are famous cases in K-theory.
Now, It is interesting that the analogus theory of commutativity preserving exterior product can be developed to the field of Lie theory. In this paper, we introduce a non abelian commutativity preserving exterior product. and the Bogomolov multiplier of Lie algebras, then we investigate their properties. Also, we compute the Bogomolov multiplier for Heisenberg Lie algebras, nilpotent Lie algebras of dimensional at most $6$ and some complex simple Lie algebras.
\section{\bf{Some notations and preliminaries}}
Let $L$ be a finite dimensional Lie algebra. The following standard notations will be used throughout the paper.
\begin{itemize}
\item{$[.,.]$                 \quad                 the Lie bracket,}
\item{$L^2=[L,L]$         \quad               the commutator subalgebra of $L$,}
\item{$H(m)$             \quad                   the Heisenberg Lie algebra of dimension $2m+1$,}
\item{$A(n)$             \quad                    the abelian Lie algebra of dimension $n$,}
\item{$\mathcal{M}(L)\cong \dfrac{R\cap{F^2}}{[R,F]}$              \quad                   the Schur multiplier of $L$, such that $L\cong \dfrac{F}{R}$.}
\\
\end{itemize}
{\bf{Exterior product 2.1.[8]}} Let $L$ be a Lie algebra and $M$ and $N$ be two ideals of $L$. the exterior product of $M$ and $N$ is defined to be the Lie algebra $M\wedge N$ generated by the symbols $m\wedge n$, where $m\in M$ and $n\in N$, subject to the following relations:
\renewcommand {\labelenumi}{(\roman{enumi})}
\begin{enumerate}
\item{$\lambda (m\wedge n) = \lambda m \wedge n = m\wedge \lambda n,$}
\item{$(m+m')\wedge n = m\wedge n + m'\wedge n,$}
\item{$m\wedge (n+n') = m\wedge n + m\wedge n',$}
\item{$[m,m']\wedge n = m\wedge [m',n] - m'\wedge [m,n],$}
\item{$m\wedge [n,n'] = [n',m]\wedge n - [n,m]\wedge n',$}
\item{$[(m\wedge n),(m'\wedge n')] = -[n,m]\wedge [m',n'],$}
\item{$m\wedge n =0 \ \ \text{whenever} \ \ m=n.$}
\end{enumerate}
for all $\lambda \in F$ , $m,m'\in M$ , $n,n'\in N$.\\
\\
{\bf{Exterior pairing 2.2.[8]}}
Let L be a Lie algebra, a function $\phi: M \times N\to L$ is called an exterior pairing, if we have
\renewcommand {\labelenumi}{(\roman{enumi})}
\begin{enumerate}
\item{$h({\lambda}m, n) = h(m, {\lambda}n) = {\lambda}h(m, n),$}
\item{$h(m + m', n) = h(m, n) + h(m', n),$}
\item{$h(m, n + n') = h(m, n) + h(m, n'),$}
\item{$h([m, m'],n) = h(m, [m', n])+ h(m',[n, m]),$}
\item{$h(m, [n, n']) = h([n', m], n) + h([m, n], n'),$}
\item{$[h(m, n), h(m', n')] = h([m, n], [m', n']),$}
\item{$h(m, n) = 0 \ \ \text{whenever}\ \  m = n.$}
\end{enumerate}
for all $\lambda \in F$ , $m,m′ \in M$ and  $n, n′ \in N$.
\\
Note that the function $M\times N\to M\wedge N$ given by $(m, n)\to m\wedge n$ is the universal exterior pairing from $M\times N$.
\section{\bf{The commutativity-preserving nonabelian exterior product of Lie algebras (CP exterior product)}}
In this section, we intend to extend the results of [4,6,14,15,17,20,26] to the theory of Lie algebras.
\\
{\bf{Definition 3.1.}}
Let $K$ be a Lie algebra and $M$ and $N$ be ideals of $K$, A bilinear function $h: M\times N \to K$, is called a Lie-${\tilde{B_0}}$-pairing if for all $m,m' \in M$ and $n,n'\in N$,
\renewcommand {\labelenumi}{(\roman{enumi})}
\begin{enumerate}
\item{$h([m,m'],n)=h(m,[m',n])-h(m',[m,n]),$}
\item{$h(m,[n,n'])=h([n',m],n)-h([n,m],n'),$}
\item{$h([n,m],[m',n'])=-[h(m,n),h(m',n')],$}
\item{If $[m,n']=0$, then $h(m,n')=0$}.
\end{enumerate}
{\bf{Definition 3.2.}}
A Lie algebra Homomorphism is a linear map $H\in \text{Hom}(L, M)$ between to Lie algebras $L$ and $M$ such that it is compatible with the Lie bracket:
$$H: L\to M \ \ \ , \ \ \ H([x,y])=[H(x),H(y)]$$
For example any vector space can be made into a Lie algebra with the trivial bracket.
\\
\\
{\bf{Definition 3.3.}}
A Lie-${\tilde{B_0}}$-pairing $h : M\times N \to L$ is said to be universal, if for any Lie-${\tilde{B_0}}$-pairing $h' : M\times N \to L'$ there is a unique Lie homomorphism $\theta : L \to L'$ such that $\theta h=h'
$.
\\The following definition extend the concept of CP exterior product in [26], to the theory of Lie algebra.
\\
\\
{\bf{Definition 3.4.}}
Let $L$ be a Lie algebra and $M$ and $N$ be ideals of $L$. The CP exterior product $M\curlywedge N$ is the Lie algebra generated by symboles $m\curlywedge n$ (for $m\in M$ , $n\in N$) subject to the following relations:
\renewcommand {\labelenumi}{(\roman{enumi})}
\begin{enumerate}
\item{$\lambda (m\curlywedge n) = \lambda m \curlywedge n = m\curlywedge \lambda n,$}
\item{$(m+m')\curlywedge n = m\curlywedge n + m'\curlywedge n,$
\\
                   $m\curlywedge (n+n') = m\curlywedge n + m\curlywedge n',$
}
\item{$[m,m']\curlywedge n = m\curlywedge [m',n] - m'\curlywedge [m,n],$
\\
                   $m\curlywedge [n,n'] = [n',m]\curlywedge n - [n,m]\curlywedge n',$}
\item{$[(m\curlywedge n),(m'\curlywedge n')] = -[n,m]\curlywedge [m',n'],$}
\item{If $[m,n]=0$, then $m\curlywedge n=0$, } \\
for all $\lambda \in F$ , $m,m' \in M$ , $n,n' \in N$.
\\
\end{enumerate}
$ $
In the case $M=N=L$, we call $L\curlywedge L$ the curly exterior product of $L$.
\\
\\
{\bf{Proposition 3.5.}}
The mapping $h : M\times N \to M\curlywedge N$ given by $(m,n) \longmapsto m\curlywedge n$, is a universal  Lie-${\tilde{B_0}}$-pairing.
\begin{proof}
By definitions 3.1, 3.3 and 3.4, the proof is straightforward.
\end{proof}
{\bf{Theorem 3.6.}}
Let $L$ be a Lie algebra and $M$ and $N$ be ideals of $L$. Then we have
$$M\curlywedge N\cong \frac{M\wedge N}{{\mathcal{M}_0}(M,N)},$$
where $\mathcal{M}_0(M,N)=<m\wedge n \ \vert \ m\in M, \ n\in N, \ [m,n]=0>$.
\begin{proof}
By using definition 2.2, the function $h: M\times N \to M\curlywedge N$ given by $(m,n)\longmapsto (m\curlywedge n)$ is an exterior pairing, so it induces a homomorphism $\tilde{h}: M\wedge N \to M\curlywedge N$, given by $\tilde{h}(m\wedge n)=m\curlywedge n$, for all $m\in M$ and $n\in N$. Clearly $\mathcal{M}_0(M,N)\subseteq \ker{\tilde{h}}$. So, we have the homomorphism ${h}^*: {(M\wedge N)}/{\mathcal{M}_0(M,N)}\to M\curlywedge N$ given by $(m\wedge n) + \mathcal{M}_0(M,N)\longmapsto (m\curlywedge n)$. On the other hand, the map ${l}^*: M\curlywedge N \to {(M\wedge N)}/{\mathcal{M}_0(M,N)}$ given by $(m\curlywedge n)\longmapsto (m\wedge n)+{\mathcal{M}_0(M,N)}$ is induced by the Lie-$\tilde{B_0}$-pairing ${l}: {M}\times {N} \to {(M\wedge N)}/{{\mathcal{M}_0(M,N)}}$ given by $(m,n)\longmapsto (m\wedge n)+\mathcal{M}_0(M,N)$. Now it is easy to see that ${{h}^*}{l^*}={l^*}{{h}^*}=1$. Thus, ${l}^*$ is an isomorphism.
\end{proof}
$ $
It is known that $\kappa: M\times N \to [M,N]$ given by $(m,n)\longmapsto [m,n]$ is a crossed pairing, so, it induces a homomorphism $\tilde{\kappa}: M\wedge N \to [M,N]$, such that $\tilde{\kappa}(m\wedge n)=[m,n]$, for all $m\in M$ and $n\in N$. The kernel of $\tilde{\kappa}$ is denoted by $\mathcal{M}(M,N)$. It can be shown that $\mathcal{M}_0(M,N) \leq \mathcal{M}(M,N)$, so that there is a homomorphism \\ ${\kappa}^*: M\wedge N / \mathcal{M}_0(M,N) \to [M,N]$ given by $m\wedge n + \mathcal{M}_0(M,N)\longmapsto [m,n]$, with $\ker {{\kappa}^*}\cong \mathcal{M}(M,N)/\mathcal{M}_0(M,N)$. Similar to groups, we denote $\mathcal{M}(M,N)/\mathcal{M}_0(M,N)$ by ${\tilde{B_0}}(M,N)$, and we call it the Bogomolov multiplier of the pair of Lie algebras $(M,N)$. Therefore, we have an exact sequence
$$0\to {\tilde{B_0}}(M,N) \to M\curlywedge N \to [M,N] \to 0.$$
In the case $M=N=L$, $\mathcal{M}_0(L,L)=<l\wedge l'   \ | \ l,l'\in L\ ,\ [l,l']=0>$ and we denote it by $\mathcal{M}_0(L)$.
\\
It is known that the kernel of ${{\kappa}}: L\wedge L \to L^2$ given by $l\wedge l'\longmapsto [l,l']$ is the Schur multiplier of $L$. On the other hand $\mathcal{M}_0(L) \leq \mathcal{M}(L) = \ker {\kappa}$. So there is a homomorphism ${\tilde{{\kappa}}} : L\wedge L / \mathcal{M}_0(L) \to L^2$ given by $l\wedge l' + \mathcal{M}_0(L)\longmapsto [l,l']$ and $\ker {\tilde{{\kappa}}}\cong \mathcal{M}(L)/\mathcal{M}_0(L)$. Similar to groups, we denote $\mathcal{M}(L)/\mathcal{M}_0(L)$ by ${\tilde{B_0}}(L)$, and we call it the Bogomolov multiplier of Lie algebra $L$. So, we have an exact sequence of Lie algebra as follows,
$$0\to {\tilde{B_0}}(L) \to L\curlywedge L \to L^2 \to 0.$$
{\bf{Proposition 3.7.}}
Let $L$ be a Lie algebra and $M$, $N$ be ideals of $L$, and $K$ be an ideal of $L$ which is contained in $M\cap N$. Then there is an isomorphism
$${M}/{K}\curlywedge {N}/{K} \cong {(M\curlywedge N)}/{T},$$
where $ T = <m\curlywedge n \ |\  m\in M  ,  n\in N  ,  [m,n]\in K>$.
\begin{proof}
The map $\phi: M\times N\to {M}/{K}\curlywedge {N}/{K}$ given by $(m,n) \to (m+K)\curlywedge (n+K)$ is a well-defined Lie-$\tilde{B_0}$-pairing. Thus there is a homomorphism ${\phi}^*: M\curlywedge N\to {M}/{K}\curlywedge {N}/{K}$ with $m\curlywedge n\longmapsto (m+K)\curlywedge (n+K)$. Clearly $T\subseteq \ker{{\phi}^*}$. So, we have the homomorphism $\psi: {(M\curlywedge N)}/{T}\to {M}/{K}\curlywedge {N}/{K}$ given by $m\curlywedge (n + T)\longmapsto (m+K)\curlywedge (n+K)$. On the other hand, the map ${\varphi}^*: {M}/{K}\curlywedge {N}/{K} \to {(M\curlywedge N)}/{T}$ given by $(m+K)\curlywedge (n+K)\longmapsto (m\curlywedge n)+T$ is induced by the Lie-$\tilde{B_0}$-pairing ${\varphi}: {M}/{K}\times {N}/{K} \to {(M\curlywedge N)}/{T}$ given by $(m+K,n+K)\longmapsto (m\curlywedge n)+T$. One can check that, ${{\varphi}^*}{\psi}={\psi}{{\varphi}^*}=1$. Thus, ${\varphi}^*$ is an isomorphism. The proof is complete.
\end{proof}
Now, we give the behavior of the CP exterior product respect to a direct sum of Lie algebras.
\\
{\bf{Proposition 3.8.}} Suppose $L_1$ and $L_2$ are ideals of a Lie algebra $L$. Then
$$(L_1\oplus L_2)\curlywedge (L_1\oplus L_2)\cong L_1\curlywedge L_1\oplus L_2 \curlywedge L_2$$
\begin{proof}
The result obtained by using a similar way to that of [8].
\\
\end{proof}
\section{\bf{Hopf-type formula for the Bogomolov multiplier of Lie algebras}}
Let $L$ be a Lie algebra with a free presention $L\cong {F}/{R}$, where $F$ is a free Lie algebra and $R$ is an ideal of $F$. By the well-known Hopf formula [8], we have an isomorphism $\mathcal{M}(L)\cong (R\cap F^2)/[R,F]$. Here we intend to give the similar formula for ${\tilde{B_0}}(L)$.
\\
\\
In the following $K(F)$ denotes $\{ [x,y] \ |\  x,y\in F \}$. \\
{\bf{Proposition 4.1.}}
Let $L$ be a Lie algebra with a free presention $L\cong {F}/{R}$. Then
$${\tilde{B_0}}(L)\cong\dfrac{R\cap F^2}{<K(F)\cap R>}.$$
\begin{proof}
From [8], we have isomorphisms $L\wedge L \cong {F^2}/[R,F]$ and $L^2\cong {F^2}/{(R\cap F^2)}$. Also, the map $\tilde{\kappa} : L\wedge L \to {L^2}$ given by $x\wedge y \to [x,y]$ is an epimorphism. Thus, $\ker \tilde{\kappa}=\mathcal{M}(L)\cong (R\cap {F^2})/[R,F]$ and $\mathcal{M}_0(L)$ can be determined with the subalgebra of ${F}/{[R,F]}$ generated by all the commutators in ${F}/{[R,F]}$ that belong to $\mathcal{M}(L)$. Thus,
$$\mathcal{M}_0(L)\cong <K(\dfrac{F}{[R,F]})\cap \dfrac{R}{[R,F]}> = \dfrac{<K(F)\cap R>+[R,F]}{[R,F]} =\dfrac{<K(F)\cap R>}{[R,F]}.$$
Therefore ${\tilde{B_0}}(L) ={\mathcal{M}(L)}/{\mathcal{M}_0(L)}\cong {R\cap F^2}/{<K(F)\cap R>}$ as required.
\end{proof}
$ $
{\bf{Proposition 4.2.}}
Let $L$ be a Lie algebra and $M$ be an ideal of $L$. Then the following sequence is exact.
$${\tilde{B_0}}(L) \to {\tilde{B_0}}(\dfrac{L}{M})\to \dfrac{M}{<K(L)\cap M>} \to \dfrac{L}{L^2}\to \frac{{L}/{M}}{({L}/{M})^2}\to 0$$
\begin{proof}
Suppose $0 \xrightarrow{} R \xrightarrow{}  F \xrightarrow{\pi} L \xrightarrow{} 0$ be a free presention of $L$ and let \\ $T=\ker(F\to L/M)$. We have $M\cong T/R$. The inclusion maps $R\cap {F^2}\xrightarrow{f} T\cap {F^2},\ \ $ $T\cap {F^2}\xrightarrow{g} T,\ \ $ $T\xrightarrow{h} F$ and $F\xrightarrow{k} F$ induce the sequence of homomorphism \\
$\displaystyle{\frac{R\cap {F^2}}{<K(F)\cap R>} \xrightarrow{f^*} \frac{T\cap {F^2}}{<K(F)\cap T>} \xrightarrow{g^*} \frac{T}{<K(F)\cap T>+R} \xrightarrow{h^*} \frac{F}{R+{F^2}} \xrightarrow{k^*}}$\\ $\displaystyle{\frac{F}{T+{F^2}} \to 0}$. Observe that
$\displaystyle{\frac{T}{<K(F)\cap T>+R}\cong \frac{M}{<K(L)\cap M>}}$, $\displaystyle{\frac{F}{R+{F^2}}\cong \dfrac{L}{L^2}}$, $\displaystyle{\frac{F}{T+{F^2}}\cong \frac{{L}/{M}}{({L}/{M})^2}}$.
Now by using theorem 4.1, we have \\
$\displaystyle{{\tilde{B_0}}(L)\cong \frac{R\cap{F^2}}{<K(F)\cap R>}}$ and $\displaystyle{{\tilde{B_0}}(\dfrac{L}{M})\cong\frac{T\cap{F^2}}{<K(F)\cap T>}}$. also,\\ $\displaystyle{\text{im}{f^*}=\ker{g^*}=}$ $\displaystyle{\frac{R\cap{F^2}}{<K(F)\cap T>}},\ \ $ $\displaystyle{\text{im}{g^*}=\ker{h^*}=\frac{T\cap{F^2}}{<K(F)\cap T>+R}}$, \\ $\displaystyle{\text{im}{h^*}=\ker{k^*}=\frac{T}{R+{F^2}}}$. and $K^*$ is epimorphism. So, above sequence is exact.
\\
\end{proof}
{\bf{Proposition 4.3.}}
Let $L$ be a Lie algebra with a free presention $L\cong {F}/{R}$, and $M$  be an ideal of $L$, such that $T=\ker(F\to {L}/{M})$. Then the sequence
$$0\to \dfrac{R\ \cap <K(F)\cap T>}{<K(F)\cap R>} \to {\tilde{B_0}}(L) \to {\tilde{B_0}}(\dfrac{L}{M})\to \dfrac{M\cap L^2}{<K(L)\cap M>} \to 0,$$
 is exact.
\begin{proof}
Suppose $0 \xrightarrow{} R \xrightarrow{}  F \xrightarrow{\pi} L \xrightarrow{} 0$ be a free presention of $L$ and let \\ $T=\ker(F\to L/M)$. We have $M\cong T/R$. The inclusion maps\\ $R\ \cap <K(F)\ \cap \ T>\xrightarrow{f} R\ \cap {F^2},\ \ $ $R\ \cap {F^2}\xrightarrow{g} T\cap {F^2}$. and the map $T\cap {F^2}\xrightarrow{h} (T\cap {F^2})+R$ induce the sequence of homomorphism \\
$\displaystyle{0\to \frac{R \cap <K(F)\cap T>}{<K(F)\cap R>}\xrightarrow{f^*}\frac{R\cap {F^2}}{<K(F)\cap R>} \xrightarrow{g^*} \frac{T\cap {F^2}}{<K(F)\cap T>} \xrightarrow{h^*}}$\\ $\displaystyle{\frac{(T\cap {F^2})+R}{<K(F)\cap T>+R} \to 0}$. By using Dedekind law observe that \\
$\displaystyle{\frac{M\cap {L^2}}{<K(L) \cap M>}\cong \frac{(T\cap {F^2})+R}{<K(F) \cap T>+R}}$. Now by using theorem 4.1, we have\\
$\displaystyle{{\tilde{B_0}}(L)\cong \frac{R\cap{F^2}}{<K(F)\cap R>}},\ \ $ $\displaystyle{{\tilde{B_0}}(\dfrac{L}{M})\cong\frac{T\cap{F^2}}{<K(F)\cap T>}},\ \ $ \\ $\displaystyle{\text{im}{f^*}=\ker{g^*}=\frac{R\ \cap<K(F)\cap T>}{<K(F)\cap R>}},\ \ $ $\displaystyle{\text{im}{g^*}=\ker{h^*}=\frac{R\cap{F^2}}{<K(F)\cap T>}}$. \\ and $h^*$ is epimorphism. So, the above sequence is exact.
\\
\end{proof}
For groups, Schur multiplier is a universal object of central extensions. Recently, parallel to the classical theory of central extensions, Jezernik and Moravec in [14 ,15] developed a version of extension that preserve commutativity. They showed that the Bogomolov multiplier is also universal object parametrizing such extension for a given group. Now, we want to introduce a similar notion for Lie algebras.
\\
\\
{\bf{Definition 4.4.}}
Let $L$, $M$ and $C$ be Lie algebras. An exact sequence of Lie algebras $0 \xrightarrow{} M \xrightarrow{\chi}  C \xrightarrow{\pi} L \xrightarrow{} 0$, is called a comutativity preserving extension (CP extension) of $M$ by $L$, if commuting pairs of elements of $L$ have commuting lifts in $C$. A special type of CP extensions with the central kernel is named a central CP extension.
\\
\\
{\bf{Proposition 4.5.}}
Let $e : 0 \xrightarrow{} M \xrightarrow{\chi}  C \xrightarrow{\pi} L \xrightarrow{} 0$ be a central extension. Then $e$ is a CP extension if and only if $\chi (M) \cap K(C)=0$.
\begin{proof}
Suppose that $e$ is a CP central extension. Let $[c_1,c_2] \in \chi (M) \cap K(C)$ Then there is a commuting lift $(c'_1,c'_2)\in C\times C$ of the commuting pair $(\pi(c_1),\pi(c_2))$, such that $\pi(c'_1)=\pi(c_1)$ and $\pi(c'_2)=\pi(c_2)$. So we have $c'_1={c_1}+a$ , $c'_2={c_2}+b$ for some $a,b \in \chi (M)$. Therefore, $0=[c'_1,c'_2]=[{c_1}+a,{c_2}+b]=[c_1,c_2]$. Hence $\chi (M) \cap K(C)=0$. Conversly, suppose that $\chi (M) \cap K(C)=0$. Choose $x,y\in L$ with $[x,y]=0$. We have $x=\pi(c_1)$ and $y=\pi(c_2)$ for some $c_1 , c_2 \in C$. Therefore $\pi ([c_1,c_2])=0$, Hence $[c_1,c_2] \in \chi (M) \cap K(C)=0$ and so $[c_1,c_2]=0$. So the central extension $e$ is a CP extension.
\\
\end{proof}
{\bf{Definition 4.6.}}
We call an abelian ideal $M$ of a Lie algebra $L$, a CP Lie subalgebra of $L$ if the extension $0\to M \to L \to \dfrac{L}{M} \to 0$ is a CP extension.
\\
Also by using proposition 4.5 an abelian ideal $M$ of a Lie algebra $L$ is a CP Lie subalgebra of $L$ if $M\cap K(L)=0$.
\\
\\
Now, we obtain an explicit formula for the Bogomolov multiplier of a direct product of two Lie algebras. The following Lemma gives a free presention for ${L_1}\oplus{L_2}$, in terms of the given free presention for $L_1$ and $L_2$. it will be used in our next investigation.
\\
\\
{\bf{Lemma 4.7. [28, Lemma 2.1]}}
Let $L_1$ and $L_2$ be Lie algebras with free presentions ${F_1}/{R_1}$ and ${F_2}/{R_2}$, respectively. and let $F = F_1\ast F_2$ be the free product of $F_1$ and $F_2$. Then $0\to R \to F \to L_1 \oplus L_2 \to 0$ is a free presention for $L_1 \oplus L_2$, where $R=R_1+ R_2+[F_2,F_1]$.
\\
\\
{\bf{Proposition 4.8.}}
Let $L_1$ , $L_2$ be two Lie algebras. Then
$${\tilde{B_0}}(L_1\oplus L_2) \cong {\tilde{B_0}}(L_1) \oplus {\tilde{B_0}}(L_2).$$
\begin{proof}
By using Lemma 4.7, we have
$${\tilde{B_0}}(L_1\oplus L_2) = \dfrac{(R_1+ R_2+[F_2,F_1])\cap (F_1\ast F_2)^2}{<K(F_1\ast F_2)\cap (R_1+ R_2+[F_2,F_1])>}$$
\\
Now, let $F = F_1\ast F_2$, then the epimorphism $F\to F_1\times F_2$ induces the following epimorphism
$$\alpha : \dfrac{R\cap F^2}{<K(F)\cap R>}\to  \dfrac{R_1\cap {F_1}^2}{<K(F_1)\cap R_1>} \oplus \dfrac{R_2\cap {F_2}^2}{<K(F_2)\cap R_2>}$$
$$x+<K(F)\cap R>\longmapsto (x_1 + <K(F_1)\cap R_1>\ ,\ x_2 + <K(F_2)\cap R_2>)$$
where $x=x_1 + x_2$, such that, $x_1\in R_1\cap {F_1}^2$ and $x_2\in R_2\cap {F_2}^2$.\\
On the other hand, the map
$$\beta : \dfrac{R_1\cap {F_1}^2}{<K(F_1)\cap R_1>} \oplus \dfrac{R_2\cap {F_2}^2}{<K(F_2)\cap R_2>} \to  \dfrac{R\cap F^2}{<K(F)\cap R>}$$
given by
$$ (x_1 + <K(F_1)\cap R_1>\ ,\ x_2 + <K(F_2)\cap R_2>) \longmapsto x+<K(F)\cap R>$$
is a well-defined homomorphism. It is easy to check that $\beta$ is a right inverse to $\alpha$, so $\alpha$ is an epimorphism. Now, we show that $\alpha$ is a monomorphism. Let $x+<K(F)\cap R>\ \in \ker \alpha$, such that, $x=t_1 + t_2$. So we have $t_1\in <K(F_1)\cap R_1>$ and $t_2\in <K(F_2)\cap R_2>$. Since $t_1 , t_2 \in <R\cap K(F)>$ then $x\in <K(F)\cap R>$. Thus $\alpha$ is a monomorphism.
\\
\end{proof}
\section{\bf{Computing the Bogomolov multiplier of Heisenberg Lie algebras}}
$ $
We use the symbol $H(m)$ for the Heisenberg Lie algebras of dimenstion $2m+1$. The Heisenberg Lie algebra $L$ is a Lie algebra such that $L^2=Z(L)$ and $\dim L^2=1$. Such Lie algebras are odd dimensional with basis $v_1,\ldots ,v_{2m}, v$ and the only non–zero multiplication between basis elements is $[v_{2i-1}, v_{2i}] = -[v_{2i}, v_{2i-1}] = v$ for $i = 1, 2,\ldots ,m$.
\\
{\bf{Theorem 5.1.}}
${\tilde{B_0}}(H(1))=0$.
\begin{proof}
Since $H(1)\wedge H(1) = <{v_1}\wedge {v_2} , {v_1}\wedge v , {v_2}\wedge v>$, an element $w\in \mathcal{M}(H(1))\leq H(1)\wedge H(1)$ can be written as $w={\alpha_1}({v_1}\wedge {v_2})+{\alpha_2}({v_1}\wedge v)+{\alpha_3}({v_2}\wedge v)$, for ${\alpha_1 ,\alpha_2 , \alpha_3} \in F$. Now, considering $\tilde{\kappa} : H(1)\wedge H(1) \to {H(1)}^2$ with $\ker {\tilde{\kappa}}=\mathcal{M}(H(1))$, we have $\tilde{\kappa} (w)=0$, and hence ${\alpha_1}[{v_1},{v_2}]+{\alpha_2}[{v_1},v]+{\alpha_3}[{v_2},v]=0$. On the other hand, $[{v_1},v]=[{v_2},v]=0$, $[{\alpha_1} {v_1},{v_2}]={\alpha_1}[{v_1},{v_2}]={\alpha_1} v=0$. Hence ${v_1}\wedge v , {v_2}\wedge v , {\alpha_1}({v_1}\wedge {v_2}) \in \mathcal{M}_0(H(1))$. Thus $w\in \mathcal{M}_0(H(1))$, and so $\mathcal{M}(H(1))\subseteq \mathcal{M}_0(H(1))$. Therefore ${\tilde{B_0}}(H(1))=0$.
$ $
\end{proof}
{\bf{Theorem 5.2.}}
${\tilde{B_0}}(H(m))=0$, for all $m\geq 2$.
\begin{proof}
We know that
$$H(m) = <v_1,v_2,\ldots v_{2m},v \ |\  [v_{2i-1},v_{2i}]=-[v_{2i},v_{2i-1}]=v   ,   i=1\ldots m>.$$
so, we can see that \\
$H(m)\wedge H(m) = <v_1\wedge v_2 , v_1\wedge v_3 , \ldots , v_1\wedge v_{2m} , v_2\wedge v_3 , v_2\wedge v_4 , \ldots , v_2 \wedge v_{2m} , \ldots $ \\
$ v_{2m-1}\wedge v_{2m}, v_1\wedge v , \ldots , v_{2m}\wedge v>$.
\\
Also for all $1\leq i\leq m$, we have\\
$v_i\wedge v = v_i \wedge [v_{2i-1},v_{2i}]=[v_{2i},v_i]\wedge v_{2i-1}-[v_{2i-1},v_i]\wedge v_{2i}=0$. Thus, \\ $H(m)\wedge H(m) = <v_1\wedge v_2 , v_1\wedge v_3 , \ldots , v_1\wedge v_{2m} , v_2\wedge v_3 , \ldots , v_2\wedge v_{2m} , \ldots , v_{2m-1}\wedge v_{2m}>$. \\
Now for all $w\in \mathcal{M}(H(m)) \leq H(m)\wedge H(m)$, there exists \\ $\alpha_1, \ldots \alpha _{2m^2-2m} , \beta_1 , \ldots , \beta _m \in F$, such that, $w=\alpha_1(v_1 \wedge v_3)+\alpha_2(v_1 \wedge v_4)+\ldots +\alpha_{2m^2-2m}(v_{2m-2}\wedge v_{2m})+\beta_1(v_1\wedge v_2)$
$+\beta_2(v_3 \wedge v_4)+\ldots +\beta_m(v_{2m-1}\wedge v_{2m})$. \\ Let $\tilde{\kappa} : H(m)\wedge H(m) \to {H(m)}^2$ is given by $x\wedge y \to [x,y]$. Since $\tilde{\kappa} (w)=0$, we have $\alpha_1[v_1 , v_3]+\alpha_2[v_1 , v_4]+\ldots +\alpha_{2m^2-2m}[v_{2m-2} , v_{2m}]+\beta_1[v_1 , v_2]+\beta_2[v_3 , v_4]+\ldots +\beta_m[v_{2m-1} , v_{2m}]=0$. So, $(\beta_1 + \beta_2 +\ldots + \beta_m)v=0$. Hence, $w=\alpha_1(v_1 \wedge v_3)+\alpha_2(v_1 \wedge v_4)+\ldots +\alpha_{2m^2-2m}(v_{2m-2}\wedge v_{2m})$
$+\beta_1(v_1\wedge v_2 - v_3\wedge v_4)+\beta_2(v_3 \wedge v_4 - v_5\wedge v_6)+\ldots +\beta_{m-1}(v_{2m-3}\wedge v_{2m-2}-v_{2m-1}\wedge v_{2m})$. \\
On the other hand, $[v_1,v_3]=[v_1,v_4]=\ldots =[v_{2m-2},v_{2m}]=0$, Thus \\
$v_1\wedge v_3 , v_1\wedge v_4 , \ldots v_{2m-2}\wedge v_{2m}\in M_0(H(m))$. We can see that, $[v_1 + v_4 , v_2 + v_3]=0$, So, $(v_1 + v_4)\wedge( v_2 + v_3)\in \mathcal{M_0}(H(m))$. Hence, $v_1\wedge v_2 + v_1\wedge v_3 +v_4\wedge v_2 +v_4\wedge v_3 \in \mathcal{M}_0(H(m))$. Thus, $(v_1\wedge v_2) - (v_3\wedge v_4) \in \mathcal{M}_0(H(m))$.\\
By a same way, we have,
$$((v_3\wedge v_4)-(v_5\wedge v_6)), \ldots ,((v_{2m-3}\wedge v_{2m-2})-(v_{2m-1}\wedge v_{2m}))\in \mathcal{M}_0(H(m)).$$ Therefore $w\in \mathcal{M}_0(H(m))$ and so $\mathcal{M}(H(m))\subseteq \mathcal{M}_0(H(m))$. Hence ${\tilde{B_0}}(H(m))=0$ as required.
\\
\end{proof}
{\bf{Lemma 5.3.}}
Let $L$ be an $n$-dimensional Lie algebra with $\dim L^2=1$. Then ${\tilde{B_0}}(L)=0$.
\begin{proof}
By Lemma 3.3 in [27], $L\cong H(m)\oplus A(n-2m-1)$ for some $m$. Now using Theorem 5.2 and Proposition 4.8, we have
$${\tilde{B_0}}(L)\cong {\tilde{B_0}}(H(m)\oplus A(n-2m-1))\cong {\tilde{B_0}}(H(m))\oplus {\tilde{B_0}}(A(n-2m-1)).$$
Since ${\tilde{B_0}}(H(m))={\tilde{B_0}}(A(n-2m-1))=0$, the result follows.
\\
\end{proof}
\section{\bf{Computing The Bogomolov multiplier of nilpotent Lie algebras of dimensional at most $6$}}
$ $
This section is devoted to obtain the Bogomolov multiplier for the nilpotent Lie algebras of dimension at most $6$. We need the classification of these Lie algebras in [7,10]. The following results are obtained by using notations and terminology in [1,6,14,16].
\\
{\bf{Theorem 6.1.}}
Let $L$ be a nilpotent Lie algebra of dimension at most $2$, Then ${\tilde{B_0}}(L)=0$.
\begin{proof}
Since $L$ is abelian, its Bogomolov multiplier is trivial.
\end{proof}
From [10], there are two nilpotent Lie algebras of dimension $3$, the abelian one, which we denote by $L_{3,1}$ and $L_{3,2}\cong H(1)$ with basis ${v,v_1,v_2}$ and nonzero Lie bracket $[v_1,v_2]=v$.
\\
\\
{\bf{Theorem 6.2.}}
Let $L$ be a nilpotent Lie algebra of dimension $3$. Then ${\tilde{B_0}}(L)=0$.
\begin{proof}
$L_{3,1}$ is abelian Lie algebra. Thus ${\tilde{B_0}}(L_{3,1})=0$. Now since $L_{3,2}\cong H(1)$. the result is obtained by using theorem 5.1.
\\
\end{proof}
From [10], there are three nilpotent Lie algebras of dimensional $4$, which are isomorphic to $L_{4,1},L_{4,2},L_{4,3}$ and  $L_{4,k}\cong L_{3,k}\oplus I   ,   k=1,2$ (where $I$ is $1$-dimensional abelian ideal). $L_{4,3}$ has the basis ${x_1,x_2,x_3,x_4}$, by non zero brackets $[x_1,x_2]=x_3$, $[x_1,x_3]=x_4$.
\\
\\
{\bf{Theorem 6.3.}}
Let $L$ be a nilpotent Lie algebra of dimension $4$. Then ${\tilde{B_0}}(L)=0$.
\begin{proof}
Using Proposition 4.8 and theorem 6.2 we have
$${\tilde{B_0}}(L_{4,k})\cong {\tilde{B_0}}(L_{3,k})\oplus {\tilde{B_0}}(I)=0, \ \text{for}\  k=1,2.$$
Let new $L\cong L_{4,3}=<x_1,x_2,x_3,x_4\  |\  [x_1,x_2]=x_3 , [x_1,x_3]=x_4>$, we have\\ $x_2\wedge x_4=x_3\wedge x_4=0$. So, $L_{4,3}\wedge L_{4,3} = <x_1\wedge x_2 , x_1\wedge x_3 , x_1\wedge x_4 , x_2\wedge x_3>$.
\\
Hence, for all $w\in \mathcal{M}(L_{4,3})\leq L_{4,3}\wedge L_{4,3}$ there exists $\alpha_1,\alpha_2,\alpha_3 ,\alpha_4\in F$, such that $w=\alpha_1(x_1\wedge x_2)+\alpha_2(x_1\wedge x_3)+\alpha_3(x_1\wedge x_4)+\alpha_4(x_2\wedge x_3)$. Now, considering \\ $\tilde{\kappa} : L_{4,3}\wedge L_{4,3} \to {L_{4,3}}^2$ given by $x\wedge y \to [x,y]$. Since $\tilde{\kappa} (w)=0$, we have $\alpha_1[x_1,x_2]+\alpha_2[x_1,x_3]+\alpha_3[x_1,x_4]+\alpha_4[x_2,x_3]=0$. So, $\alpha_1 x_3 + \alpha_2 x_4=0$. On the other hand, $[x_1,x_4]=[x_2,x_3]=[x_2,x_4]=[x_3,x_4]=0$ , $[\alpha_1 x_1,x_2]={\alpha_1}[x_1,x_2]=\alpha_1 x_3=0$ and $[\alpha_2 x_1,x_3]={\alpha_2}[x_1,x_3]=\alpha_2 x_4=0$. Hence
$(x_1\wedge x_4) , (x_2\wedge x_3) , \alpha_1(x_1\wedge x_2) ,$  $\alpha_2(x_1\wedge x_3)\in \mathcal{M}_0(L_{4,3})$. So, $\mathcal{M}(L_{4,3}) \subseteq \mathcal{M}_0(L_{4,3})$. Thus ${\tilde{B_0}}(L_{4,3})=0$.
\\
\end{proof}
From [10], The $5$-dimensional Lie algebras are $L_{5,k}\cong L_{4,k} \oplus I$, for $k=1,2,3$. where $I$ is a $1$-dimensional abelian ideal and the following Lie algebras
\begin{itemize}
\item{$L_{5,4}\cong <x_1,...,x_5 \ | \ [x_1,x_2]=[x_3,x_4]=x_5>,$}
\item{$L_{5,5}\cong <x_1,...,x_5 \ | \ [x_1,x_2]=x_3 , [x_1,x_3]=[x_2,x_4]=x_5>,$}
\item{$L_{5,6}\cong <x_1,...,x_5 \ | \ [x_1,x_2]=x_3 , [x_1,x_3]=x_4 , [x_1,x_4]=[x_2,x_3]=x_5>,$}
\item{$L_{5,7}\cong <x_1,...,x_5 \ | \ [x_1,x_2]=x_3 , [x_1,x_3]=x_4 , [x_1,x_4]=x_5>,$}
\item{$L_{5,8}\cong <x_1,...,x_5 \ | \ [x_1,x_2]=x_4 , [x_1,x_3]=x_5>,$}
\item{$L_{5,9}\cong <x_1,...,x_5 \ | \ [x_1,x_2]=x_3 , [x_1,x_3]=x_4 , [x_2,x_3]=x_5>.$}
\\
\end{itemize}
{\bf{Theorem 6.4.}}
Let $L$ be a nilpotent Lie algebra of dimension $5$. Then ${\tilde{B_0}}(L)\neq 0$ if and only if $L\cong L_{5,6}$.
\begin{proof}
By using theorem 6.3, Proposition 4.8, one can check that \\
${\tilde{B_0}}(L_{5,1})={\tilde{B_0}}(L_{5,2})={\tilde{B_0}}(L_{5,3})={\tilde{B_0}}(L_{5,4})={\tilde{B_0}}(L_{5,7})={\tilde{B_0}}(L_{5,8})={\tilde{B_0}}(L_{5,9})=0$.
Now let $L\cong L_{5,5}$, we can see that
$$L_{5,5}\wedge L_{5,5} = <x_1\wedge x_2 , x_1\wedge x_3 , x_1\wedge x_4 , x_2\wedge x_3 , x_2\wedge x_4 , x_3\wedge x_4>.$$
Hence, for all $w\in \mathcal{M}(L_{5,5})\leq L_{5,5}\wedge L_{5,5}$, there exists $\alpha_1,\alpha_2,\ldots , \alpha_6 \in F$, such that
$w=\alpha_1(x_1\wedge x_2)+\alpha_2(x_1\wedge x_3)+\alpha_3(x_1\wedge x_4)+\alpha_4(x_2\wedge x_3)+\alpha_5(x_2\wedge x_4)+\alpha_6(x_3\wedge x_4)$. \\
Since $\tilde{\kappa} (w)=0$, we have $\alpha_1[x_1,x_2]+\alpha_2[x_1,x_3]+\alpha_3[x_1,x_4]+\alpha_4[x_2,x_3]+\alpha_5[x_2,x_4]+\alpha_6[x_3,x_4]=0$. Thus, $\alpha_1 x_3 + (\alpha_2+\alpha_5) x_5=0$. Therefore $\alpha_1 x_3=(\alpha_2+\alpha_5) x_5=0$. Hence, $w=\alpha_1(x_1\wedge x_2)+\alpha_2((x_1\wedge x_3)-(x_2\wedge x_4))+\alpha_3(x_1\wedge x_4)+\alpha_4(x_2\wedge x_3)$. \\ On the other hand $[x_1,x_4]=[x_2,x_3]=[\alpha_1 x_1,x_2]=\alpha_1[x_1,x_2]=\alpha_1 x_3=0$. So, $(x_1\wedge x_4) , (x_2\wedge x_3) , {{\alpha}_1}(x_1\wedge x_2)\in \mathcal{M}_0(L_{5,5})$. We can see that, $[x_1+x_2+x_3,x_1+x_2+x_4]=0$ and so $(x_1+x_2+x_3)\wedge (x_1+x_2+x_4)\in \mathcal{M}_0(L_{5,5})$, Hence $(x_1\wedge x_4)+(x_2\wedge x_4)+(x_3\wedge x_1)+(x_3\wedge x_2)+(x_3\wedge x_4)\in \mathcal{M}_0(L_{5,5})$ and $(x_1\wedge x_3)-(x_2\wedge x_4) \in \mathcal{M}_0(L_{5,5})$. Therefore $\mathcal{M}(L_{5,5})\subseteq \mathcal{M}_0(L_{5,5})$, and hence ${\tilde{B_0}}(L_{5,5})=0$. Let $L\cong L_{5,6}$, we can see that $L_{5,6}\wedge L_{5,6} = <x_1\wedge x_2 , x_1\wedge x_3 , x_1\wedge x_4 , x_1\wedge x_5 , x_2\wedge x_3 , x_2\wedge x_5>$. Hence, for all $w\in \mathcal{M}(L_{5,6})\leq L_{5,6}\wedge L_{5,6}$ there exists $\alpha_1,\alpha_2,\ldots , \alpha_6 \in F$, such that
$w=\alpha_1(x_1\wedge x_2)+\alpha_2(x_1\wedge x_3)+\alpha_3(x_1\wedge x_4)+\alpha_4(x_1\wedge x_5)+\alpha_5(x_2\wedge x_3)+\alpha_6(x_2\wedge x_5)$. \\
Since $\tilde{\kappa} (w)=0$, we have $\alpha_1[x_1,x_2]+\alpha_2[x_1,x_3]+\alpha_3[x_1,x_4]+\alpha_4[x_1,x_5]+\alpha_5[x_2,x_3]+\alpha_6[x_2,x_5]=0$. Thus, $\alpha_1 x_3 +\alpha_2 x_4 +(\alpha_3+\alpha_5) x_5=0$. Therefore $\alpha_1 x_3=\alpha_2 x_4=(\alpha_3+\alpha_5) x_5=0$. On the other hand, $[\alpha_1 x_1,x_2]=\alpha_1[x_1,x_2]=\alpha_1 x_3=0$
and $[\alpha_2 x_1,x_3]=\alpha_2[x_1,x_3]=\alpha_2 x_4=0$.
So, $\alpha_1(x_1\wedge x_2) , \alpha_2(x_1\wedge x_3) , (x_1\wedge x_5) , (x_2\wedge x_5) \in \mathcal{M}_0(L_{5,6})$. Thus, $w=\alpha_3(x_1\wedge x_4-x_2\wedge x_3)+\tilde{w}$, where $\tilde{w}\in \mathcal{M}_0(L_{5,6})$. Let $A$ be a generating set for $\mathcal{M}_0(L_{5,6})$, then $\mathcal{M}(L_{5,6})=<A,\ (x_1\wedge x_4-x_2\wedge x_3)>$. Hence, $\dim \tilde{B_0}(L_{5,6})=1$. So ${\tilde{B_0}}(L_{5,6})\cong A(1)$.
\end{proof}
From [10], The $6$-dimensional Lie algebras are $L_{6,k}\cong L_{5,k} \oplus I$, for $k=1,\ldots ,9$. where $I$ is a $1$-dimensional abelian ideal and following Lie algebras
\begin{itemize}
\item{$L_{6,10}\cong <x_1,...,x_6 \ | \ [x_1,x_2]=x_3 , [x_1,x_3]=[x_4,x_5]=x_6>,$}
\item{$L_{6,11}\cong <x_1,...,x_6\ | \ [x_1,x_2]=x_3 , [x_1,x_3]=x_4 , [x_1,x_4]=[x_2,x_3]=[x_2,x_5]=x_6>$,}
\item{$L_{6,12}\cong <x_1,...,x_6 \ | \ [x_1,x_2]=x_3 , [x_1,x_3]=x_4 , [x_1,x_4]=[x_2,x_5]=x_6>,$}
\item{$L_{6,13}\cong <x_1,...,x_6 \ | \ [x_1,x_2]=x_3 , [x_1,x_3]=[x_2,x_4]=x_5 , [x_1,x_5]=[x_3,x_4]=x_6>,$}
\item{$L_{6,14}\cong <x_1,...,x_6 \ | \ [x_1,x_2]=x_3 , [x_1,x_3]=x_4 , [x_1,x_4]=[x_2,x_3]=x_5 , $
\\
$[x_2,x_5]=x_6 , [x_3,x_4]=-x_6>,$}
\item{$L_{6,15}\cong <x_1,...,x_6 \ | \ [x_1,x_2]=x_3 , [x_1,x_3]=x_4 , [x_1,x_4]=[x_2,x_3]=x_5 ,$
\\
$[x_1,x_5]=[x_2,x_4]=x_6>,$}
\item{$L_{6,16}\cong <x_1,...,x_6 \ | \ [x_1,x_2]=x_3 , [x_1,x_3]=x_4 , [x_1,x_4]=x_5 , [x_2,x_5]=x_6 , [x_3,x_4]=-x_6>,$}
\item{$L_{6,17}\cong <x_1,...,x_6 \ | \ [x_1,x_2]=x_3 , [x_1,x_3]=x_4 , [x_1,x_4]=x_5 , [x_1,x_5]=[x_2,x_3]=x_6>,$}
\item{$L_{6,18}\cong <x_1,...,x_6 \ | \ [x_1,x_2]=x_3 , [x_1,x_3]=x_4 , [x_1,x_4]=x_5 , [x_1,x_5]=x_6>,$}
\item{$L_{6,19}(\epsilon)\cong <x_1,...,x_6 \ | \ [x_1,x_2]=x_4 , [x_1,x_3]=x_5 , [x_2,x_4]=x_6 , [x_3,x_5]={\epsilon}x_6>$, ($L_{6,19}(\epsilon)\cong L_{6,19}(\delta)$ if and only if there is an $\alpha \in {F^*}$ such that $\delta={{\alpha}^2}{\epsilon}$).}
\item{$L_{6,20}\cong <x_1,...,x_6 \ | \ [x_1,x_2]=x_4 , [x_1,x_3]=x_5 , [x_1,x_5]=[x_2,x_4]=x_6>,$}
\item{$L_{6,21}(\epsilon)\cong <x_1,...,x_6 \ | \ [x_1,x_2]=x_3 , [x_1,x_3]=x_4 , [x_2,x_3]=x_5 , [x_1,x_4]=x_6 , [x_2,x_5]={\epsilon}x_6>$, ($L_{6,21}(\epsilon)\cong L_{6,21}(\delta)$ if and only if there is an $\alpha \in {F^*}$ such that $\delta={{\alpha}^2}{\epsilon}$).}
\item{$L_{6,22}(\epsilon)\cong <x_1,...,x_6 \ | \ [x_1,x_2]=x_5 , [x_1,x_3]=x_6 , [x_2,x_4]={\epsilon}x_6 , [x_3,x_4]=x_5>$, ($L_{6,22}(\epsilon)\cong L_{6,22}(\delta)$ if and only if there is an $\alpha \in {F^*}$ such that $\delta={{\alpha}^2}{\epsilon},$).}
\item{$L_{6,23}\cong <x_1,...,x_6 \ | \ [x_1,x_2]=x_3 , [x_1,x_3]=[x_2,x_4]=x_5 , [x_1,x_4]=x_6>,$}
\item{$L_{6,24}(\epsilon)\cong <x_1,...,x_6 \ | \ [x_1,x_2]=x_3 , [x_1,x_3]=[x_2,x_4]=x_5 , [x_1,x_4]={\epsilon}x_6 , [x_2,x_3]=x_6>$, ($L_{6,24}(\epsilon)\cong L_{6,24}(\delta)$ if and only if there is an $\alpha \in {F^*}$ such that $\delta={{\alpha}^2}{\epsilon}$).}
\item{$L_{6,25}\cong <x_1,...,x_6 \ | \ [x_1,x_2]=x_3 , [x_1,x_3]=x_5 , [x_1,x_4]=x_6>,$}
\item{$L_{6,26}\cong <x_1,...,x_6 \ | \ [x_1,x_2]=x_4 , [x_1,x_3]=x_5 , [x_2,x_3]=x_6>.$}
\end{itemize}
$ $
{\bf{Theorem 6.5.}}
Let $L$ be a nilpotent Lie algebra of dimension $6$. Then ${\tilde{B_0}}(L)\neq 0$ if and only if $L$ is isomorphic to one of Lie algebra $L_{6,6}$, $L_{6,13}$, $L_{6,14}$, $L_{6,15}$, $L_{6,19}(e)$; $(e\geq 1)$, $L_{6,20}$, $L_{6,21}(1)$, $L_{6,22}(0)$, $L_{6,23}$, $L_{6,24}(e)$; $(e\geq 0)$.
\begin{proof}
By using a Similar method involving in Theorem 6.4 the results follow.
\end{proof}
One of the important results for the Schur multiplier was presented by Moneyhun in [25]. He showed that for a Lie algebra $L$ of dimension $n$, $\dim \mathcal{M}(L) = n(n-1)/2-t(L)$. for some $t(L)\geq 0$. His results suggestes an interesting problem: Can we classify Lie algebras of dimension $n$ by $t(L)$? the answer to this question can found for $t(L)=1,\ldots, 8$ in [2,12,13,25]. On the other hand, from [27], we have an upper bound for the dimension of the Schur multiplier of a non-abelian nilpotent Lie algebra as follows $\dim \mathcal{M}(L) = n(n-1)(n-2)/2 +1-s(L)$, for some $s(L)\geq 0$. Hence by the same motivation we have the analogues question for clasification of $L$ according to $s(L)$. It seems classifying nilpotent Lie algebras by $s(L)$ helps the classification of Lie algebras in term of $t(L)$. (See for instance [27]). Now, according to this classification, we will investigate Bogomolov multiplier for some Lie algebras.
\\
\\
{\bf{Theorem 6.6.}}
Let $L$ be an $n$-dimensional nilpotent Lie algebra with $s(L)=1$, Then ${\tilde{B_0}}(L)=0$.
\begin{proof}
Since $s(L)=1$, by Theorem 3.9 in [27], $L\cong L_{5,4}$. So, ${\tilde{B_0}}(L)=0$.
\end{proof}
$ $
{\bf{Theorem 6.7.}}
Let $L$ be an $n$-dimensional nilpotent Lie algebra and $t(L)\leq 6$, then  ${\tilde{B_0}}(L)=0$.
\begin{proof}
By Theorem 3.10 in [27] and also Proposition 4.7, ${\tilde{B_0}}(L)=0$.
\end{proof}
$ $
{\bf{Theorem 6.8.}}
Let $L$ be an $n$-dimensional nilpotent Lie algebra with $s(L)=2$ and $dimL^2=2$. Then ${\tilde{B_0}}(L)=0$.
\begin{proof}
By Theorem 4.3 in [27], $L\cong L_{4,3}$ or $L\cong L_{5,4}\oplus A(1)$, thus ${\tilde{B_0}}(L)=0$.
\\
\end{proof}
\section{\bf{The Bogomolov multipliers of some complex simple Lie algebras}}
$ $
A simple group is a group with no nontrivial proper normal subgroup. and the classification of finite simple groups is a major milestone in the history of mathematics. On the other hand with the help of the Jordan-Holder theorem, a finite group can be written as a certain combination of simple groups. Also, in contrast to the classification of finite simple groups, the classification of simple Lie groups is simplified by using the manifold structure. In particular every Lie group has an dependent Lie algebra, and in this regard, some authors have also gained some results. for example, Bosshardt showed that a Lie group is simple if and only if its Lie algebra is simple. (see [5,30,31] for more information). Theories of groups and Lie algebras are structurally similar, and many concepts related to groups, there are analogously defined concepts for Lie algebras. Eventually, this subject reduces the problem of finding simple Lie groups to classifying simple Lie algebras. and in this section we obtain the Bogomolov multiplier of some complex simple Lie algebras.
\\
\\
{\bf{Definition 7.1.[11]}}
A Lie algebra $L$ is simple if it has no ideals other than $0$ and $L$, and it is not abelian.
\\
\\
{\bf{Definition 7.2.[11]}}
A Lie algebra $L$ is called semisimple if the only commutative ideal of $L$ is $0$. for example $0$-dimensional Lie algebra, the special linear Lie algebra, the odd-dimensional special orthogonal Lie algebra, the symplectic Lie algebra and the even-dimensional special orthogonal Lie algebra for $(n>1)$ are semisimple.
\\
\\
Complex simple Lie algebras have been completely classified by Cartan [5]. They classified into four infinite classes with five exceptional Lie algebras.
\\
\\
{\bf{Theorem 7.3.[11]}}
Every simple Lie algebra over $\Bbb{C}$ is isomorphic to precisely one of the following Lie algebras:
\\
1. $Sl(n+1; \Bbb{C})$   ,   $n\geq 1$,
\\
2. $So(2n+1 ; \Bbb{C})$  ,  $n\geq 2$,
\\
3. $Sp(n ; \Bbb{C})$  ,  $n\geq 3$,
\\
4. $So(2n ; \Bbb{C})$   ,  $n\geq 4$,
\\
5. The exceptional Lie algebras $G_2$ , $F_4$ , $E_6$ , $E_7$ and $E_8$.
\\
\\
Knapp in [19] showed the five exceptional Lie algebras $G_2, F_4, E_6, E_7, E_8$ have dimension $14, 52, 78, 133$ and $248$, respectively.
\\
\\
In the following $E_{ij}$ denotes the matrix with $1$ at the intersection of the $i$-th row and the $j$-th coloumn and $0$ every where else. The Lie bracket of $E_{ij}$ and $E_{kl}$ is given by $[E_{ij},E_{kl}]={E_{ij}}{E_{kl}}-{E_{kl}}{E_{ij}}={{\delta}_{jk}}{E_{il}}-{{\delta}_{il}}{E_{kj}}$.
\\
\\
{\bf{Theorem 7.4.}} Let $L$ be a special linear Lie algebras $sl(n+1; \Bbb{C})$. Then $\tilde{B_0}(L)=0$.
\begin{proof}
From [19], $sl(n+1; \Bbb{C})$ has the basis $D_{{i}{i+1}}$, $E_{ij}$ such that $D_{ij}=E_{ii}-E_{jj}$. So, for $j\neq i =1\ldots n$, we have
\\
$sl(n+1; \Bbb{C})=<D_{i{i+1}} , E_{ij}\ |\ [D_{i{i+1}},D_{{i+1}{i+2}}]= D_{i{i+2}} , [D_{i{i+1}},E_{ij}]=2E_{ij} \ ;\ $\\ $ j=i+1 \  ,\
[D_{i{i+1}},E_{ij}]=E_{ij}\ ;\ j\neq i+1>$
$\mod \mathcal{M}_0(sl(n+1; \Bbb{C}))$.
\\
We can see that
$sl(n+1; \Bbb{C})\wedge sl(n+1; \Bbb{C})=<D_{i{i+1}}\wedge D_{{i+1}{i+2}} , D_{i{i+1}}\wedge E_{ij}>$  \\ mod $\mathcal{M}_0(sl(n+1; \Bbb{C}))$. Now, for all
$w\in \mathcal{M}(sl(n+1; \Bbb{C})) \leq sl(n+1; \Bbb{C})\wedge   sl(n+1; \Bbb{C})$, there exists ${\alpha}_{i} , {\beta}_{ij} \in \Bbb{C} , i,j =1...n+1$ and $\tilde{w}\in \mathcal{M}_0(sl(n+1; \Bbb{C}))$ such that
$$w={\sum^{n}_{i=1}}{{\alpha}_i}(D_{i{i+1}}\wedge D_{{i+1}{i+2}})+{\sum^{n}_{i,j=1}}{{\beta}_{ij}}(D_{i{i+1}}\wedge E_{ij}) + {\tilde{w}}.$$
Since $\tilde{\kappa} (w)=0$, we have ${\sum^{n}_{i=1}}{{\alpha}_i}[D_{i{i+1}}\wedge D_{{i+1}{i+2}}]+{\sum^{n}_{i,j=1}}{{\beta}_{ij}}[D_{i{i+1}}\wedge E_{ij}]=0$. Now, if $j=i+1$, then ${\sum^{n}_{i=1}}{{\alpha}_i}{D_{i{i+2}}}+{\sum^{n}_{i=1}}2{{\beta}_{i}}{E_{i{i+1}}}={\sum^{n}_{i=1}}{{\alpha}_i}(E_{ii}-E_{{i+2}{i+2}})+2{\sum^{n}_{i=1}}{{\beta}_i}E_{i{i+1}}=0$.
So, for all $i=1...n$, ${{\alpha}_i}={{\beta}_i}=0$. If $j\neq i+1$, then ${\sum^{n}_{i=1}}{{\alpha}_i}D_{i{i+2}}+{\sum^{n}_{i,j=1 , i<j}}{{\beta}_{i,j}}E_{ij}=0$. So, for all $i,j$, ${{\alpha}_i}={{\beta}-{ij}}=0$. Hence, \\$w\in \mathcal{M}_0(sl(n+1; \Bbb{C}))$ and $\mathcal{M}(sl(n+1; \Bbb{C}))\subseteq \mathcal{M}_0(sl(n+1; \Bbb{C}))$. Therefore $\tilde{B}_0(sl(n+1; \Bbb{C}))=0$.
\end{proof}
$ $
{\bf{Theorem 7.5.}} Let $L$ be one of the odd-dimensional orthogonal Lie algebras $so(2n+1; \Bbb{C})$. Then $\tilde{B_0}(L)=0$.
\begin{proof}
From [19], $so(2n+1; \Bbb{C})$ has the basis $H_i$, ${K_i}^{\pm}$, ${L_{ij}}^{\pm}$, ${M_{ij}}^{\pm}$ such that
$ $
$D_{ij}=E_{ij}-E_{ji}$,     $1\leq i\neq j\leq 2n+1$
$ $
$H_i := \sqrt{-1}D_{{2i-1}{2i}}$,    $i=1,...,n$
$ $
${K_i}^{\pm} := D_{{2i-1}{2n+1}} \pm \sqrt{-1}D_{{2i}{2n+1}}$,    $i=1,...,n$
$ $
${L_{ij}}^{\pm} := (D_{{2i-1}{2j-1}}-D_{{2i}{2j}}) \pm \sqrt{-1}(D_{{2i-1}{2j}}+D_{{2i}{2j-1}})$,    $1\leq i<j\leq n$
$ $
${M_{ij}}^{\pm} := (D_{{2i-1}{2j}}-D_{{2i}{2j-1}}) \pm \sqrt{-1}(D_{{2i-1}{2j-1}}+D_{{2i}{2j}})$,    $1\leq i<j \leq n$.
$ $
Also, we have
$ $
$[H_i,{K_{i}}^{\pm}]=\sqrt{-1}D_{{2n+1}{2i}}\pm D_{{2i-1}{2n+1}}$
$ $
$[H_i,{L_{ij}}^{\pm}]=-\sqrt{-1}D_{{2i-1}{2j}}\pm D_{{2i-1}{2j-1}}$
$ $
$[H_i,{M_{ij}}^{\pm}]=-\sqrt{-1}D_{{2i-1}{2j-1}}\pm D_{{2i-1}{2j}}$
$ $
$[{K_{i}}^{\pm} , {L_{ij}}^{\pm}]=[{K_{i}}^{\pm} , {M_{ij}}^{\pm}]=[{L_{ij}}^{\pm} , {M_{ij}}^{\pm}]=0$. So in$\mod \mathcal{M}_0(so(2n+1; \Bbb{C}))$, we can see that \\
$so(2n+1; \Bbb{C})=<H_i , {K_i}^{\pm} , {L_{ij}}^{\pm} ,  {M_{ij}}^{\pm}\ |\ [H_i,{K_{i}}^{\pm}],[H_i,{L_{ij}}^{\pm}],[H_i,{M_{ij}}^{\pm}]>$
\\
 and $so(2n+1; \Bbb{C})\wedge so(2n+1; \Bbb{C}) = <H_i \wedge {K_{i}}^{\pm} , H_i \wedge {L_{ij}}^{\pm} , H_i \wedge {M_{ij}}^{\pm}>$
\\
$=<D_{{2i-1}{2i}} \wedge D_{{2i}{2n+1}} , D_{{2i-1}{2i}} \wedge D_{{2i}{2j}}>$  mod $\mathcal{M}_0(so(2n+1; \Bbb{C}))$. Now for all \\ $w\in \mathcal{M}(so(2n+1; \Bbb{C})) \leq so(2n+1; \Bbb{C})\wedge so(2n+1; \Bbb{C})$, there exists ${\alpha}_1 , {\alpha}_2 \in \Bbb{C}$ and $\tilde{w} \in \mathcal{M}_0(so(2n+1; \Bbb{C}))$, such that $w={{\alpha}_1}(D_{{2i-1}{2i}} \wedge D_{{2i}{2n+1}})+{{\alpha}_2}(D_{{2i-1}{2i}} \wedge D_{{2i}{2j}}) +{\tilde{w}}$.
\\
Since $\tilde{\kappa} (w)=0$, we have \\ ${{\alpha}_1}[D_{{2i-1}{2i}} , D_{{2i}{2n+1}}]+{{\alpha}_2}[D_{{2i-1}{2i}} , D_{{2i}{2j}}]={{\alpha}_1}{D_{{2i-1}{2n+1}}}+{{\alpha}_2}{D_{{2i-1}{2j}}}$ \\
$={{\alpha}_1}(E_{{2i-1}{2n+1}}-E_{{2n+1}{2i-1}})+{{\alpha}_2}(E_{{2i-1}{2j}}-E_{{2j}{2i-1}})=0$. \\
Thus, ${\alpha}_1={\alpha}_2=0$. Hence, $\mathcal{M}(so(2n+1; \Bbb{C}))\subseteq \mathcal{M}_0(so(2n+1; \Bbb{C}))$ and \\ $\tilde{B}_0(so(2n+1; \Bbb{C}))=0$.
\end{proof}
$ $
{\bf{Theorem 7.6.}} Let $L$ be an even-dimensional orthogonal Lie algebras $so(2n; \Bbb{C})$, Then $\tilde{B}_0(L)=0$.
\begin{proof}
From [19], $so(2n; \Bbb{C})$ has the basis $H_i$, ${L_{ij}}^{\pm}$, ${M_{ij}}^{\pm}$. such that
$ $
$D_{ij}=E_{ij}-E_{ji}$,      $1\leq i\neq j\leq 2n+1$
$ $
$H_i := \sqrt{-1}D_{{2i-1}{2i}}$,    $i=1,...,n$
$ $
${L_{ij}}^{\pm} := (D_{{2i-1}{2j-1}}-D_{{2i}{2j}}) \pm \sqrt{-1}(D_{{2i-1}{2j}}+D_{{2i}{2j-1}})$,    $1\leq i<j\leq n$
$ $
${M_{ij}}^{\pm} := (D_{{2i-1}{2j}}-D_{{2i}{2j-1}}) \pm \sqrt{-1}(D_{{2i-1}{2j-1}}+D_{{2i}{2j}})$,    $1\leq i<j \leq n$.
$ $
Also, we have
$ $
$[H_i,L_{ij}]=-\sqrt{-1}D_{{2i-1}{2j}}\pm D_{{2i-1}{2j-1}}$
$ $
$[H_i,M_{ij}]=-\sqrt{-1}D_{{2i-1}{2j-1}}\pm D_{{2i-1}{2j}}$
$ $
$[L_{ij},M_{ij}]=0$
$ $
Thus, $so(2n; \Bbb{C})=<H_i , {L_{ij}}^{\pm} ,  {M_{ij}}^{\pm}\ |\ [H_i,L_{ij}],[H_i,M_{ij}]> \mod \mathcal{M}_0(so(2n; \Bbb{C}))$. \\ We can see that $so(2n; \Bbb{C})\wedge so(2n; \Bbb{C}) = < H_i \wedge {L_{ij}}^{\pm} , H_i \wedge {M_{ij}}^{\pm}>$ mod $\mathcal{M}_0(so(2n; \Bbb{C}))$. Now for all $w\in \mathcal{M}(so(2n; \Bbb{C})) \leq so(2n; \Bbb{C})\wedge so(2n; \Bbb{C})$, there exists ${\alpha}_1 , {\alpha}_2 \in \Bbb{C}$ and $\tilde{w} \in \mathcal{M}_0(so(2n; \Bbb{C}))$, such that $w={{\alpha}_1}(H_i \wedge {L_{ij}}^{\pm})+{{\alpha}_2}(H_i \wedge {M_{ij}}^{\pm})+{\tilde{w}}$. Since $\tilde{\kappa} (w)=0$, we have ${{\alpha}_1}[H_i , {L_{ij}}^{\pm}]+{{\alpha}_2}[H_i , {M_{ij}}^{\pm}]=
{{\alpha}_1}(-\sqrt{-1}D_{{2i-1}{2j}}\pm D_{{2i-1}{2j-1}})+$
\\
${{\alpha}_2}(-\sqrt{-1}D_{{2i-1}{2j-1}}\pm D_{{2i-1}{2j}})=0$.
Thus, ${\alpha}_1={\alpha}_2=0$.
Hence, \\ $\mathcal{M}(so(2n; \Bbb{C}))\subseteq \mathcal{M}_0(so(2n; \Bbb{C}))$ and so $\tilde{B}_0(so(2n; \Bbb{C}))=0$.
\end{proof}
$ $
{\bf{Theorem 7.7.}} Let $L$ be a sympletic Lie algebras $sp(n; \Bbb{C})$. Then $\tilde{B}_0(L)=0$.
\begin{proof}
From [19], $sp(n; \Bbb{C})$ has the basis $H_i$, $X_{ij}$, $Y_{ij}$, $Z_{ij}$, $U_i$, $V_i$. such that
$ $
$H_i=E_{ii}-E_{{n+i}{n+i}}$,    $1\leq i \leq n$
$ $
$X_{ij} := {E_{ij}}-{E_{{n+j}{n+i}}}$,    $1\leq i\neq j\leq n$
$ $
 $Y_{ij} := {E_{i{n+j}}}+{E_{j{n+i}}}$,     $1\leq i<j\leq n$
$ $
$Z_{ij} := {E_{{n+i}j}}+{E_{{n+j}i}}$,    $1\leq i<j\leq n$
$ $
$U_i := {E_{i{n+i}}}$,     $1\leq i \leq n$
$ $
$V_i := {E_{{n+i}i}}$,    $1\leq i \leq n$.
Since
$ $
$[X_{ij},Y_{ij}]=2{E_{i{n+i}}},\ [X_{ij},Z_{ij}]=-2E_{{n+j}j},\ [X_{ij},V_i]=-E_{{n+i}j}-E_{{n+j}i},\ $ $[Y_{ij},Z_{ij}]=E_{ii}+E_{jj},\ [Y_{ij},V_i]=-E_{{n+i}{n+j}}+E_{ji},\ [Z_{ij},U_i]=-E_{ij}+E_{{n+j}{n+i}},\ [U_i,V_i]=E_{ii},\ [X_{ij},U_i]=[Y_{ij},U_i]=[Z_{ij},V_i]=0,\ [H_i,X_{ij}]=-E_{{n+j}{n+i}},\ [H_i,Y_{ij}]=E_{i{n+j}}-E_{j{n+i}},\ [H_i,Z_{ij}]=-E_{{n+i}j}-E_{{n+j}i},\ [H_i,U_i]=2E_{i{n+i}},\ [H_i,V_i]=-2E_{{n+i}i}$,
we have
\\
$sp(n; \Bbb{C})=<H_i,X_{ij},Y_{ij},Z_{ij},U_{i},V_{i} \ |\ [X_{ij},Y_{ij}],[X_{ij},Z_{ij}],[X_{ij},V_{i}],[Y_{ij},Z_{ij}],[Y_{ij},V_{i}],$\\$[Z_{ij},U_{i}],[U_{i},V_{i}],[H_i,X_{ij}],[H_i,Y_{ij}],[H_i,Z_{ij}],[H_i,U_i],[H_i,V_i]>$\\$\ \mod  \mathcal{M}_0(sp(n; \Bbb{C}))$. One can see that
\\
$sp(n; \Bbb{C})\wedge sp(n; \Bbb{C}) \equiv <{X_{ij}}\wedge {Y_{ij}},{X_{ij}}\wedge {Z_{ij}},{X_{ij}}\wedge {V_i},{Y_{ij}}\wedge {Z_{ij}},{Y_{ij}}\wedge {V_i},{Z_{ij}}\wedge {U_i},$
\\
${U_i}\wedge {V_i},{H_i}\wedge {X_{ij}},{H_i}\wedge {Y_{ij}},{H_i}\wedge {Z_{ij}},{H_i}\wedge {U_i},{H_i}\wedge {V_i}> \mod  \mathcal{M}_0(sp(n; \Bbb{C}))$.
\\
By using a similar method involving in previous theorems, the result follow.
\\
\end{proof}

\end{document}